\documentclass{article}

\input{amssym.def}
\input{amssym.tex}

\usepackage{amsthm}
\usepackage{epsfig}
\usepackage{color}

\title{Heegaard genus and Property $\tau$ for hyperbolic 3-manifolds}
\author{D. D. Long\thanks{supported in part by the NSF},~A. Lubotzky\thanks{supported in part by the NSF and the US-Israel
Bi-National Science Foundation}~ \&~ A. W. Reid\thanks{supported in part 
by the NSF}}

\def\ker{\rm{ker}}

\def\inf{\rm{inf}}
\def\mod{\rm{mod}}
\def\PSL{\rm{PSL}}
\def\SL{\rm{SL}}
\def\GL{\rm{GL}}
\def\qed{ $\sqcup\!\!\!\!\sqcap$}

\newtheorem{theorem}{Theorem}[section]
\newtheorem{lemma}[theorem]{Lemma}

\newtheorem{proposition}[theorem]{Proposition}

\begin{document}
\maketitle

%
%
%
%
\section{Introduction}
~~~Let $M$ be a finite volume hyperbolic $3$-manifold and ${\cal L} =
\{M_i \}$ some family of finite sheeted regular covers of $M$.  We say that
${\cal L}$ is {\em co-final} if $\bigcap_i \pi_1(M_i) = \{ 1 \}$,
where, as usual, the $\pi_1(M_i)$ are all to be regarded as subgroups
of $\pi_1(M)$.  The {\em infimal Heegaard gradient} of $M$ with
respect to the family ${\cal L}$ is defined as:
$$\inf_i {\chi_-^h(M_i)\over {[\pi_1(M):\pi_1(M_i)]}},$$
where
$\chi_{-}^h(M_i)$ denotes the minimal value for
the negative of the Euler characteristic of a Heegaard surface in
$M_i$.  

In \cite{La1} Lackenby showed that if $\pi_1(M)$ is an arithmetic
lattice in $\PSL(2,{\bf C})$ then $M$ has a co-final family of covers
(namely, those arising from congruence subgroups) with positive
infimal Heegaard gradient.  The main point of the current note is to
show that the same applies for every finite volume hyperbolic
$3$-manifold:

\begin{theorem} 
\label{heegaardgrad}
Let $M$ be a finite volume hyperbolic $3$-manifold.  Then
M has a co-final family of finite sheeted covers for which the infimal
Heegaard gradient is positive.  \end{theorem}

It is interesting to recall that if $M$ has a finite sheeted cover which
fibers over a circle (and a well known conjecture due to Thurston, asserts
that every hyperbolic $M$ has such a covering) then $M$ has a (co-final)
family of finite sheeted covers whose infimal Heegaard gradient is zero.

Theorem \ref{heegaardgrad} is a consequence of Theorem \ref{newmain} below.  
To state this theorem, we need some preliminary definitions.
Let $\Gamma$ be a group generated by some finite symmetric set $S$ and 
let ${\cal L} = \{N_i \}$ be a family of finite index normal subgroups 
of $\Gamma$.
Then the group $\Gamma$ is said to have Property $\tau$ with respect to
${\cal L}$ if the family of Cayley graphs $X(\Gamma/N_i,S)$ forms a family
of expanders (see \cite{Lu} Chapter 4
for various equivalent forms of Property $\tau$).

In \cite{La1} Lackenby showed that for $M$ a finite volume hyperbolic
$3$-manifold and ${\cal L} =\{M_i\}$ a family of finite sheeted covers
(but not necessarily a cofinal family), then
if $\pi_1(M)$ has
Property $\tau$ with respect to $\pi_1(M_i)$, $M_i\in {\cal L}$,
then the infimal Heegaard gradient of $M$
with respect to $\cal L$ is positive.  Thus Theorem \ref{heegaardgrad}
follows immediately from our next result, 
which can be viewed as providing a first step towards
a generalization of Clozel's result \cite{Cl} to non-arithmetic
lattices.

\begin{theorem} 
\label{newmain}
Let $\Gamma$ be a finitely generated non-elementary Kleinian group.
Then $\Gamma$ has a co-final family of finite index normal subgroups ${\cal L}
= \{N_i\}$ with respect to which $\Gamma$ has Property $\tau$.
\end{theorem}

In \S 2, we review Property $\tau$ and prove a lifting result about
Property $\tau$ from smaller subgroups (see Proposition \ref{promote_tau})
and in \S 3, we prove 
Theorem \ref{newmain}. The
proof of Theorem \ref{newmain} 
uses a recent result of Bourgain and Gamburd \cite{BG} on expanding
properies of the finite groups $\SL(2,p)$ when the generating set makes
the Cayley graph have large girth.
This result can then be combined with Lemma \ref{margulis}, which generalises 
a result of Margulis  \cite{Ma} and shows that in our context, we can
ensure that the Cayley graph does have large girth.
It is perhaps worth emphasizing that our main result has a purely topological
conclusion,
but the methods of \cite{BG} are those of
additive combinatorics, and
in particular recent work of Tao and Vu \cite{TV},  work of Helfgott
\cite{He}, and sum-product estimates
in finite fields (\cite{BGK} and \cite{BKT}).

The key tool in Theorem \ref{newmain}  is proving a somewhat
stronger result, namely Proposition \ref{algmain}, which shows that
co-final families of the required type exist for every non-virtually soluble
subgroup of $\SL(2,k)$, where $k$ is an arbitrary number field. Theorem 
\ref{newmain} will now follow: if the Kleinian group $\Gamma$ has finite co-covolume,
then Mostow-Prasad rigidity implies that $\Gamma$ admits a discrete faithful 
representation into $\SL(2,{\bf C})$ where the entries lie in some number field.
The general case of a finitely generated Kleinian group follows by applying 
results of Brooks,  Scott and Thurston.

Finally, in \S 4 we discuss some possible generalizations of Theorem 
\ref{newmain} to finding nested co-final families and to
more general linear groups.

\section{Promoting $\tau$}

Let $\Gamma$ be a group, $S$ a finite symmetric set of generators of
$\Gamma$, ${\cal L} = \{N_i\}$ a family of finite index normal
subgroups of $\Gamma$, and $X(\Gamma/N_i,S)$ the quotient Cayley graphs.
Recall the definition that $X(\Gamma/N_i,S)$ form a family of expanders 
(see \cite{Lu}).\\[\baselineskip]
\noindent{\bf Definition:}~{\em A finite k-regular graph $X=X(V,E)$ 
with a set $V$ of n vertices is called
an $(n,k,c)$-expander if for every subset $A\subset V$:
$$|\partial A| \geq c\biggl(1 - {|A|\over n}\biggr )|A|,$$
where $\partial A = \{v\in V:d(v,A)=1\}$ and $d$ the distance function on the
graph.  

A family of k-regular graphs ($k$ fixed) is called an expander family 
if there is a $c>0$ such that all of the 
graphs are $(n,k,c)$-expanders.}\\[\baselineskip]
In this notation the family of Cayley graphs $X(\Gamma/N_i,S)$ are
$([\Gamma:N_i],|S|,c)$-expanders for some $c>0$.

There are various methods of lifting Property $\tau$ from smaller
subgroups of a larger group to the larger group (see \cite{BS} and
\cite{Lu1}).  These have been used in the context of arithmetic
groups. We now provide another method (inspired by Example E on p.52
of \cite{Lu}) which applies to both arithmetic and non-arithmetic
groups.

\begin{proposition}
\label{promote_tau}
Let $\Gamma$ be a finitely generated group, and ${\cal L} = \{N_i\}$ a family of finite
index normal subgroups of $\Gamma$. Suppose that $H<\Gamma$ (not
necessarily of finite index) is finitely generated, and assume that $H$ surjects onto the
finite quotients $\Gamma/N_i$ for all but a finite number of $i$. Then
if $H$ has Property $\tau$ with respect to the family
$\{H\cap N_i\}$, $\Gamma$
has Property $\tau$ with respect to $\cal L$.\end{proposition}

\noindent{\bf Proof:}~Fix a symmetric generating set $S_H$ for $H$ and 
extend this to a generating set $S$ for $\Gamma$. 
By assumption, the quotient Cayley graphs 
$X(H/H\cap N_i,S_H)$ form a family of expanders. Now
$H$ surjects onto the
finite quotient $\Gamma/N_i$ for all but a finite number of $i$ and so we can
ignore this finite number for our considerations 
since throwing away a finite
number of graphs will not affect the fact that a family forms a set of
expanders. 
It follows that  $X(\Gamma/N_i, S)$ is also a family of expanders, since for
every subset $A$ of $\Gamma/N_i$, the computation of the ratio
$|\partial A|/ \biggl((1 - {|A|\over n_i})|A|\biggr )$ 
(where $n_i = [\Gamma:N_i] = [H:H\cap N_i]$)
can only be increased by the
addition of the extra edges coming from enlarging the generating set.\qed

\section{Proof of Theorem \ref{newmain}}

The key result in the proving of Theorem \ref{newmain} is our next proposition.

\begin{proposition}
\label{algmain}
Let $k$ be a number field and let $\Gamma$ be a finitely generated subgroup 
of $\SL(2,k)$ which is not virtually soluble. Then 
$\Gamma$ has a co-final family of finite index normal subgroups
${\cal L}=\{N_i\}$ with respect to which 
$\Gamma$ has Property $\tau$.\end{proposition}

\noindent{\bf Proof:}~Since $\Gamma$ is finitely generated,
we can assume that $\Gamma < \SL(2,A)$ where $A$ is a ring of
S-integers in a number field $k$, of degree $n$ say, over $\bf Q$.
We fix some notation.  Let $\cal P$ be a prime ideal of $A$ with
residue class field $\bf F$, and let 
$$\pi_{\cal P} : \SL(2,A) \longrightarrow \SL(2,{\bf F}),$$
be the reduction homomorphism.

Since $\Gamma$ is not virtually soluble, it contains a non-abelian
free subgroup $F = <a^{\pm 1}, b^{\pm 1}>$.
By the Cebotarev density theorem, we can find
infinitely many rational primes $p$ which split completely in
$k$, and 
it follows from Strong Approximation (see \cite{LR}
for an elementary argument in the case of $\SL(2)$) that for all but
finitely many of the rational primes $p$ that split completely in $k$, $F$
surjects $\SL(2,p)$ under the homomorphisms $\pi_{\cal P}$. Also, note
that $\ker~\pi_{\cal P} \cap \Gamma$ (resp. $\ker~\pi_{\cal P} \cap
F$) form a co-final family of normal subgroups of finite index in
$\Gamma$ (resp.  $F$).

We need the following lemma (cf. \cite{Ma}).  Recall that the {\em
girth} of a finite graph $X$ is the length of the shortest
non-trivial closed path in $X$.

\begin{lemma}
\label{margulis}
There is a constant $C = C(a,b)$ so that the girth of the Cayley graph of
$\SL(2,p)$ with respect
to the generating set $\{\pi_{\cal P}(a^{\pm 1}) , \pi_{\cal P}(b^{\pm1 })\}$ 
is at least $C\log(p)$.
\end{lemma}

This lemma will be proved below. Assuming \ref{margulis}, the proof of 
Proposition \ref{algmain} is
completed by the following result of 
Bourgain and Gamburd (see Theorem 3 of \cite{BG}).

\begin{theorem}
\label{BG}
Suppose that for each $p$, $S_p$ is some symmetric generating set for
$\SL(2,p)$ such that the girth of the Cayley graph $X(\SL(2,p), S_p$)
at least $C\log(p)$ (where $C$ is independent of $p$). Then
$X(\SL(2,p), S_p$) form a family of expanders.\end{theorem}
%
%
%
Thus to complete the proof, applying Lemma \ref{margulis} to the free
subgroup $F$ of $\Gamma$ (with generating set as before), we can apply
Theorem \ref{BG} to $\SL(2,p)$ with these generating sets and the
result now follows from Proposition \ref{promote_tau} and our previous
discussion.\qed\\[\baselineskip]
\noindent{\bf Proof of Lemma \ref{margulis}:}\\[\baselineskip]
Denote the ring of integers of $k$ by $R_k$.  If $\alpha \in R_k$,  
define:
$$\mu(\alpha) = \max\{|\alpha'|~:~\alpha'~\hbox{is a Galois 
                        conjugate of}~\alpha\}.$$
Here $|*|$ denotes the complex absolute value. Note that
since $\alpha$ is an algebraic
integer, then $\mu(\alpha) \geq 1$ with equality if and only if $\alpha$
is a root of unity.
 
It follows easily from the definition that $\mu(\alpha + \beta) \leq
\mu(\alpha) + \mu(\beta)$
and $\mu(\alpha \cdot \beta) \leq \mu(\alpha) \cdot \mu(\beta)$, since for
example, in computing
the maximum for $|\alpha' + \beta'|$ one clearly cannot do better than
maximize the two terms
of  $|\alpha' |+ |\beta'|$.
 
Given a matrix $t \in \SL(2,A)$, we may write $t$ as $1/\alpha \cdot t^*$ 
where $t^* \in M(2,R_k)$ and $\alpha \in R_k$. 
Take $M$ to be the biggest value of $\mu$
taken over  the entries of the matrices $a^*$, $(a^{-1})^*$, $b^*$
and $(b^{-1})^*$, together with the four denominators of those matrices.\\[\baselineskip]
\noindent {\bf Claim:}~{\em If $w$ is a word of length $r$ in matrices 
taken from $a^*$,
$(a^{-1})^*$, $b^*$ and $(b^{-1})^*$, then the entries
of $w$ cannot have their $\mu$ value being larger 
than $(2M)^r$.}\\[\baselineskip]
The proof  of the claim is by induction on $r$: Consider $X\cdot a^*$ for
example, where $X$ is a word
of length $r-1$. The entries of the product  $X\cdot a^*$  have the form
$x_1a_1 + x_2a_2$, where $x_i $ is an entry in $X$ etc. Then by the above
remarks
$\mu(x_1a_1 + x_2a_2) \leq \mu(x_1a_1) + \mu(x_2a_2) \leq
\mu(x_1)\mu(a_1)+ \mu(x_2)\mu(a_2)$
and by induction this is at most $(2M)^{r-1}M + (2M)^{r-1}M  = (2M)^r$.
This completes the proof of the claim.\\[\baselineskip]
Now suppose that $p >>  0$ is a rational prime that splits completely 
in $k$ and $ \cal P $ is a $k$-prime 
dividing $p$. Let $r$ denote the girth, and let
$w(a,b)\in F$ be a
reduced word of length r which projects to a cycle of length $r$
under $\pi_{\cal P}$. Clearing denominators
in the congruence $w(a,b)= id~\mod~{\cal P}$, we obtain a congruence between 
$R_k$-integral matrices
$$ w(a^*, (a^{-1})^*, b^*, (b^{-1})^*) = Z \cdot id~~\mod~{\cal P}.$$
~~~By the claim, the entries on the lefthand side have their $\mu$ values
being bounded
above by $(2M)^r$, and $Z$ is a product of $r$ integers with $\mu$ value
at most $M$, so $\mu(Z) \leq M^r$.

Now the integral matrix $w(a^*, (a^{-1})^*, b^*, (b^{-1})^*) - Z \cdot id$
is not identically zero since $a$ and $b$ generate a free group of
rank two. Let  $\beta$ be one
of its nonzero entries. The above remarks show that
$\mu(\beta)$ is bounded above by $(2M)^r + M^r < (3M)^r$ say.
Notice that $\beta \in {\cal P}$, and the $k/{\bf Q}$-norm of $\beta$ is a
non-zero integer which is divisible  by $p$, since we take the product of
all the conjugates of $\beta$ which
lie in conjugates of $\cal P$. Recalling that  $n$ is the degree of $k$ over the rationals,
it  follows from the definition of $\mu$ that
this integer is bounded above by $\mu(\beta)^n \leq ((3M)^r)^n$.
 
We therefore deduce that in order to be divisible by $p$, 
$r$ must be large enough so that $(3M)^{rn} > p$; that is to say 
$r \geq C\log(p)$ with $C=1/(n\log(3M))$ as required.\qed\\[\baselineskip]
\noindent{\bf Proof of Theorem \ref{newmain}:}\\[\baselineskip]
Without any loss of generality, we may assume that $\Gamma$ is torsion free. In the
case that $\Gamma$ has finite co-volume, then as remarked in \S 1, it 
follows from local rigidity that $\Gamma$ (or more precisely a lift to
$\SL(2, {\bf C})$) can be conjugated into $\SL(2,k)$ for some number field $k$
(indeed $k$ can be chosen to be a quadratic extension of the trace-field, see
\cite{MR} Corollary 3.2.4), and Proposition \ref{algmain} applies.

Now assume that $\Gamma$ has infinite co-volume. Since $\Gamma$ is
non-elementary it is not virtually soluble, and is either
geometrically finite or geometrically infinite. In the former case we
can apply a result of Brooks \cite{Br1} that produces a
quasi-conformal conjugate $\Gamma'$ of $\Gamma$ that is a subgroup of
a Kleinian group of finite co-volume. Hence, as in the previous
paragraph, this implies $\Gamma'<\SL(2,k)$ for some number field $k$.
In the case when $\Gamma$ is geometrically infinite, as pointed out in
\cite{And} for example, it is a consequence of the Scott Core Theorem
and Thurston's Hyperbolization Theorem for Haken manifolds that there
exists a group $\Gamma'$ isomorphic to $\Gamma$ with $\Gamma'$
geometrically finite. We now argue as before.\qed

\section{Final Remarks}

\noindent{\bf 1.}~Although Theorem \ref{newmain} provides a co-final family,
it does not provide a family which is co-final and {\em nested}.  
In the proof of Theorem \ref{newmain} we used \cite{BG} which does not
provide a nested family. Recent work of Bourgain, Gamburd and Sarnak \cite{BGS}
extends some of \cite{BG} to products of primes, and may yet eventually lead to
a co-final nested family in groups $\Gamma$ as above.  

However, even in the absence of this result, we can
produce a rich class of non-arithmetic Kleinian groups that do contain a
co-final nested family. 

\begin{proposition}
\label{nested}
Let $\Gamma$ be a Kleinian group of finite co-volume that contains an
arithmetic Fuchsian subgroup.  Then $\Gamma$ contains a
co-final nested family
${\cal L} =\{N_i\}$ of normal subgroups of finite index such that $\Gamma$
has Property $\tau$ with respect to $\cal L$.\end{proposition}
\noindent{\bf Proof:}~We begin with a preliminary remark. 
Let $F$ be an arithmetic Fuchsian group. As is pointed out in 
\cite{Lu} pp 51 - 52 for example, 
it follows
from a result of Selberg and an application of
the Jacquet-Langlands correspondence, that
$F$ has Property $\tau$ with respect to the entire family 
of its congruence subgroups.  

Thus if $F$ is an arithmetic Fuchsian subgroup  of $\Gamma$, then after 
possibly discarding perhaps finitely many
prime ideals, we can form ``congruence subgroups of $\Gamma$'' obtained
via reduction homomorphisms. Furthermore,  by considering only the
primes that split completely in the invariant trace-field of $\Gamma$
(and hence also in the invariant trace-field of $F$), we can arrange
that there are infinitely many congruence quotients of $\Gamma$ for which
$F$ surjects. 

In fact, an easy argument using Strong Approximation and the Chinese
Remainder Theorem shows that
we can find a sequence of completely split primes, so that $F$ surjects
each of the reductions modulo  the descending sequence of ideals 
$I_n={\cal P}_1\ldots{\cal P}_n$. These now form  a co-final nested 
family of congruence subgroups. The result now follows from  
Proposition \ref{promote_tau} with $F = H$. \qed\\[\baselineskip]
Examples
of non-arithmetic Kleinian groups $\Gamma$ satisfying the hypothesis of
Proposition \ref{nested}
are plentiful as we now discuss.\\[\baselineskip]
\noindent{\bf Construction of Examples:}\\[\baselineskip]
\noindent{\bf Example 1:}~A thrice-punctured sphere has a unique hyperbolic
structure arising as ${\bf H}^2/\Gamma(2)$ where $\Gamma(2)$ is the principal
congruence subgroup of level $2$ in the modular group. In particular
$\Gamma(2)$ is an arithmetic Fuchsian group. Many non-arithmetic 
link complements contain an immersed (or embedded) thrice punctured sphere.
For example all the non-arithmetic hyperbolic twist 
knot complements. By \cite{HS} these are all mutually 
incommensurable.\\[\baselineskip]
\noindent{\bf Example 2:}~The examples of non-arithmetic hyperbolic 
manifolds of Gromov and Piatetski-Shapiro \cite{GPS} as hybrids of arithmetic 
ones. By construction these contain an arithmetic 
totally geodesic surface.\\[\baselineskip]
\noindent{\bf Example 3:}~One can also easily obtain closed examples by
only 3-dimensional methods using the 
construction of \cite{NR}. There, a non-compact
hyperbolic orbifold with a torus cusp is constructed, the boundary of which
consists of two totally geodesic isometric copies of a hyperbolic 2-orbifold
${\bf H}^2/\Delta$ for some hyperbolic triangle 
group $\Delta$.
Doubling this orbifold along the totally geodesic boundary
gives a two cusped hyperbolic orbifold ${\bf H}^3/\Gamma$
for which $\Gamma$ contains triangle groups. This construction can be done
where the triangle group is chosen arithmetic (see \cite{NR}). These 
are rigid, so any surgeries on the cusps will leave them as totally geodesic 
and arithmetic. Sufficiently high order Dehn surgeries 
will produce non-arithmetic hyperbolic 3-orbifolds (for example
by using Borel's result
on the discreteness of the set of co-volumes of arithmetic Kleinian
groups, see \cite{MR} Chapter 11.2.1). \\[\baselineskip]
\noindent{\bf 2.}~It seems quite possible that
Theorem \ref{newmain} and Proposition \ref{algmain} hold
for any finitely generated subgroup of $\SL(2,{\bf C})$ 
which is not virtually soluble. At present, our methods only
prove the following: 

\begin{theorem}
\label{weaktau}
Let $\Gamma$ be a finitely generated subgroup of $\SL(2,{\bf C})$ 
which is not virtually soluble. Then $\Gamma$
has an infinite family ${\cal L}=\{N_i\}$ of finite index normal subgroups, such that
$\Gamma$ has Property $\tau$ with respect to $\cal L$.\qed\end{theorem}
This result follows by noting that
Proposition \ref{algmain} implies that we are done unless
$\Gamma$ contains an element whose trace is transcendental; in this case
we may choose an algebraic specialization where the image group
is not virtually soluble (see \cite{LM} Proposition 2.2 for a more general
version of this specialization result).

At present our argument works only for $\SL(2)$ since
there is no analogue of the result of 
\cite{BG} yet known for groups $\SL(n)$, $n> 2$.
However, it seems reasonable to expect the following stronger conjecture to 
hold.\\[\baselineskip]
\noindent{\bf Conjecture:}~{\em Let $\Gamma$ be a finitely 
generated subgroup of $\GL(n,{\bf C})$ whose Zariski Closure is semi-simple. 
Then $\Gamma$
has a co-final (nested) family ${\cal L}=\{N_i\}$ of finite index normal subgroups,  for which 
$\Gamma$ has Property $\tau$ with respect to $\cal L$.}\\[\baselineskip]
This Conjecture would provide a far reaching generalization of 
Clozel's work \cite{Cl} mentioned in \S 1.

\bigskip
 \noindent
 Department of Mathematics,\\ University of California,\\ Santa Barbara, CA
93106, USA.

\noindent Email:~long@math.ucsb.edu\\[\baselineskip]
Department of Mathematics,\\ The Hebrew University,\\ 
Givat Ram, Jerusalem 91904, Israel.

\noindent Email:~alexlub@math.huji.ac.il\\[\baselineskip]
Department of Mathematics,\\
University of Texas\\
Austin, TX 78712, USA.

\noindent Email:~areid@math.utexas.edu\\[\baselineskip]
{\bf Abstract} We show that any finitely generated non-elementary Kleinian group has a
co-final family of finite index normal subgroups with respect to which it
has Property $\tau$.  As a consequence, any closed hyperbolic $3$-manifold
has a co-final family of finite index normal subgroups for which the
infimal Heegaard gradient is positive.

\end{document}